\documentclass{article}
\usepackage{amsmath,amssymb,amsthm,latexsym}

\theoremstyle{plain}
\newtheorem{lemma}{Lemma}[section]
\newtheorem{proposition}[lemma]{Proposition}
\newtheorem{theorem}[lemma]{Theorem}

\theoremstyle{definition}
\newtheorem{definition}[lemma]{Definition}

\theoremstyle{remark}
\newtheorem{remark}[lemma]{Remark}

\newcommand{\Q}{\mathbb{Q}}
\newcommand{\R}{\mathbb{R}}
\newcommand{\C}{\mathbb{C}}

\newcommand{\mf}[1]{\mathfrak{#1}}
\newcommand{\gn}{\mathfrak{g}}
\newcommand{\ad}{\mathop{\mathrm{ad}}}
\newcommand{\GL}{\mathop{\mathrm{GL}}}
\newcommand{\SL}{\mathop{\mathrm{SL}}}
\newcommand{\Sym}{\mathop{\mathrm{Sym}}}

\newcommand{\Tr}{\mathop{\mathrm{Tr}}}
\newcommand{\Ker}{\mathop{\mathrm{Ker}}}
\newcommand{\im}{\mathop{\mathrm{Im}}}
\newcommand{\Gal}{\mathop{\mathrm{Gal}}}
\newcommand{\PGL}{\mathop{\mathrm{PGL}}}
\newcommand{\Aut}{\mathop{\mathrm{Aut}}}
\newcommand{\compsymbol}[3]
{\vcenter to 12pt{\hsize 15pt
    \begin{center}%
      \vskip-3pt\vrule width0pt\hfill\hbox{$#1$}\hfill\vrule
width0pt\par\vskip#3\par\vrule width0pt\hfill\hbox{$#2$}\hfill\vrule width0pt%
    \end{center}%
}}
\newcommand{\iso}{\compsymbol{\sim}{\rightarrow}{-9pt}}
\newcommand{\p}{\mathbb{P}}

\title{A Lie Algebra Method for Rational Parametrization 
	of Severi-Brauer Surfaces}

\author{Willem A. de Graaf, RICAM, Linz, Austria \and 
	Michael Harrison, University of Sydney, Australia \and 
	Jana P\'\i lnikov\'a, RICAM, Linz, Austria\protect \and
	Josef Schicho, RICAM, Linz, Austria}
\begin{document}

\maketitle

\begin{abstract}
It is well-known that a Severi-Brauer surface has a rational point if
and only if it is isomorphic to the projective plane. Given a Severi-Brauer
surface, we study the problem to decide whether such an isomorphism
to the projective plane, or such a rational point, does exist; and to
construct such an isomorphism or such a point in the affirmative case.
We give an algorithm using Lie algebra techniques. The algorithm
has been implemented in Magma.
\end{abstract}

\section{Introduction}\label{se:intro}

The problem considered in this paper is to decide whether a given
surface is isomorphic to the projective plane over the rational numbers
and if so, to find an isomorphism.
It is easy to decide this over the complex numbers 
(see Section~\ref{se:antican}).
Hence we can assume that the surface is a twist of $\p^2$, also called
a Severi-Brauer surface.

The problem comes from the parametrization of the surfaces. When trying to
parametrize a surface over the rational numbers, one can reduce to several
base cases (\cite{josef00}). The Severi-Brauer surface arises 
as one of them. Therefore our problem appears as a subproblem of 
finding a parametrization of a surface with rational coefficients.

Parametrizing a Severi-Brauer surface 
(i.e.~finding an isomorphism to $\p^2$)
is equivalent to finding a rational point on the surface:
from an isomorphism one can construct all points, and in the other
direction one can construct an isomorphism when a single point is known
(\cite{manin}). This fact is not used in our paper.

There is a well-known correspondence between Severi-Brauer surfaces
and central simple algebras of degree 3 (cf.~\cite{jac96}). 
The split Severi-Brauer surfaces (i.e~those isomorphic to the projective plane)
correspond to the split central simple algebras (those isomorphic
to the full matrix algebra). 
There are a lot of classical number-theoretical results available which are
useful for deciding whether the given central simple algebra is split or not.
Here we reduce the problem to the case of cyclic algebra -- this is possible 
because of a result of Wedderburn (cf.~\cite{jac96}) -- and solve a norm
equation (\cite{pz89,fie97}).

There are known constructions
of the Severi-Brauer surface corresponding to a given central simple
algebra. However, in the other direction there are no constructions available.
We introduce an intermediate step: the Lie algebra of the given surface.
This is the Lie algebra of the group of automorphisms of the surface.

Incidentally this is also the Lie algebra of regular vector fields.
The relation between vector fields and the Lie algebra of the group of
automorphisms has been mentioned in~\cite{hm93} for the affine and local
analytic case. 

The whole algorithm has been implemented in Magma~\cite{magma}.
The most expensive step is the solution of the norm equation.

The paper is formulated for rational numbers but 
anything written generalizes to number fields. In the last step 
one needs to solve a norm equation over the number field which is also
implemented in Magma.

Most parts of the method can be extended to Severi-Brauer varieties of
any dimension. In particular we can construct the corresponding 
central simple algebra. However, the construction of the cyclic algebra
and the norm equation does not generalize.

The paper is structured as follows:
In section~\ref{se:antican} we reduce the given surface to a Del Pezzo surface
of degree 9. The isomorphism, if it exists, is then a linear projective map
to the anticanonical embedding of $\p^2$.
In section~\ref{se:lie} we reduce the problem to finding an isomorphism
of a given Lie algebra and the Lie algebra $\mf{sl}_3(\Q)$.
In section~\ref{se:sln} we reduce the problem to finding an isomorphism
of a given central simple algebra and the full matrix algebra $M_3(\Q)$.
In section~\ref{se:assoc} we describe how to solve this problem by
reducing to a norm equation.

The paper is carefully written such that the sections can be read 
in an arbitrary order. They are independent of each other.

The third author was supported by SFB Grant F1303 of the Austrian FWF. 
The authors are grateful to the Magma group of the University of Sydney 
for a stimulating research environment during a short period of the work.
We also acknowledge some helpful discussions with G\'abor Ivanyos on the
material of section~\ref{se:sln} and Peter Mayr on the material of
section~\ref{se:assoc}.

\section{The Anticanonical Embedding}\label{se:antican}

We start with a projective surface (i.e.~a variety of dimension 2) $S$ 
over the rational numbers.
%We assume that we have given a projective embedding $S\subset{\mathbb P}^r$
%for some $r$, i.e. we know a set of generators for the homogeneous
%ideal $I(S)\subset{\mathbb Q}[x_0,\dots,x_r]$. 
Our goal is to decide whether $S$ is isomorphic to ${\mathbb P}^2$;
and if yes, we want to construct an isomorphism.

A first obvious necessary condition for $S$ to be isomorphic to ${\mathbb P}^2$
is that $S$ be nonsingular. This can be checked easily by the
Jacobian criterion (see \cite{hartsh77}, Theorem~I.5.3).
So we assume from now on that $S$ is nonsingular.

For any nonsingular variety $X$, the {\em anticanonical bundle} $\cal A$
is the determinant bundle of the tangent bundle $TX$
(see \cite{manin}).
In the surface case, this is just the antisymmetric tensor bundle
of two times the tangent bundle. It is clear that $\cal A$
is always a line bundle.

For any line bundle $E$ over a projective variety $X$, the vector space
$\Gamma(X,E)$ of sections is finitely generated over the ground
field (which is $\Q$ in our case). If $\dim(\Gamma(X,E))=:n>0$,
then there exists an {\em associated rational map}
$a_E:X\to {\mathbb P}^{n-1}$, which maps the point $p$ to the
ratio $(s_1(p){:}\dots {:}s_n(p))$, where $\{s_1,\dots,s_n\}$
is a basis of $\Gamma(X,E)$. It is defined on the
complement of the subset of $X$ where all global sections vanish.

The definition of the associated map depends on
the choice of the basis, but in a transparent way: if we choose a
different basis, we get a projectively equivalent map.

The anticanonical map is the rational map associated to $\cal A$.
For the explicit construction of the
anticanonical map, we refer to \cite{deei02}.

\begin{proposition} \label{prop:s0}
The vector space of sections of the anticanonical bundle of ${\mathbb P}^2$
has dimension~10. The anticanonical map is an embedding.
It is given by
\begin{equation}\label{SBVpoint}
  (s{:}t{:}u) \mapsto (x_0{:}\dots{:}x_9) =
  (s^3{:}t^3{:}u^3{:}s^2t{:}t^2u{:}u^2s{:}st^2{:}tu^2{:}us^2{:}stu) . 
\end{equation}
% \begin{equation}\label{SBVpoint}
%   (s{:}t{:}u) \mapsto 
% 	(x_0{:}\dots{:}x_9)=(s^3{:}s^2t{:}s^2u{:}\dots{:}u^3) . 
% \end{equation}
The image is a surface $S_0$ of degree~9, whose ideal is generated by
27 quadric polynomials.
\end{proposition}

\begin{proof}
The anticanonical bundle is isomorphic to ${\cal O}(3)$
(see \cite{hartsh77}, Example~II.8.20). 
The sections correspond to cubic forms, of
which the above is a basis. 

The ideal of the image is generated by the kernel of the linear evaluation map
from the quadric forms in $x_0,\dots,x_9$ to the forms of degree~6
in $s,t,u$. This kernel has dimension $55-28=27$.
\end{proof}

Recall that a projective surface $S\subset{\mathbb P}^n$ is called
a {\em Del Pezzo surface} iff it is anticanonically embedded
(see \cite{manin}). It is well-known that $3\le n\le 9$ in this
case, and that the degree of the surface is then also $n$.

\begin{theorem} \label{thm:S}
If $S$ is isomorphic to ${\mathbb P}^2$, then 
the anticanonical map $a_{\cal A}$ is an embedding,
and the image is projectively equivalent to $S_0$.
\end{theorem}

\begin{proof}
An isomorphism $f:{\mathbb P}^2\to S$ induces
an isomorphism of the anticanonical bundles and a vector space
isomorphism of global sections. 
\end{proof}

If $X,Y$ are varieties defined over $\Q$, then we say that $X$ is a 
{\em twist} of $Y$ iff $X\otimes\C$ is $\C$-isomorphic to $Y\otimes\C$, 
(see \cite{poonen97}).
Moreover, if $X\subset{\mathbb P}^n$ and $Y\subset{\mathbb P}^n$ 
are projective varieties, $n>0$, then we say that $X$ is a 
{\em projective twist} of $Y$ iff $X$ and $Y$ are projectively
equivalent. 

If $S$ is a twist of ${\mathbb P}^2$, then its anticanonical
embedding $a_{\cal A}(S)$ is a projective twist of $S_0$,
by the complex version of Theorem~\ref{thm:S}.
The following theorem makes it possible to decide whether this statement holds
or not. 

\begin{theorem}
$S$ is a twist of ${\mathbb P}^2$ iff $\dim(\Gamma(S,{\cal A}))=10$
and the anticanonical map is an embedding.
\end{theorem}

\begin{proof}
It suffices to prove in the complex case that $S$ is isomorphic to 
${\mathbb P}^2$ iff $\dim(\Gamma(S,{\cal A}))=10$ and
the anticanonical map is an embedding.
The ``only if'' direction follows from Proposition~\ref{prop:s0}.
Conversely, the statement that the anticanonical map is an embedding
implies that $S$ is a Del Pezzo surface (see \cite{fujita90}).
By the classification of Del Pezzo surfaces (see \cite{manin}),
the equality $\dim(\Gamma(S,{\cal A}))=10$ implies that $S$
is a Severi-Brauer surface.
\end{proof}

Assume that $\dim(\Gamma(S,{\cal A}))=10$. Then we know that
$S$ is a twist of ${\mathbb P}^2$, and its anticanonical embedding 
$S_1:=a_{\cal A}(S)$ is a projective twist of $S_0$. Moreover,
$S$ is equivalent to ${\mathbb P}^2$ iff 
$S_1$ is projectively equivalent to $S_0$.
Clearly, every projective transformation $S_0\to S_1$ can
be composed with the inverse of the anticanonical embedding
to give an isomorphism ${\mathbb P}^2\to S$.
Hence, we have reduced the original problem to deciding whether
a given Del Pezzo surface of degree~9 is projectively equivalent
to $S_0$, and to compute a projective transformation in the
affirmative case.

%xxx
%In terms of polynomials, the $3$-uple embedding of ${\mathbb P}^2$
%is given by
%In fact, the ideal is generated by the kernel of the linear evaluation map
%from the quadric forms in $x_0,\dots,x_9$ to the forms of degree~6
%in $s,t,u$. This kernel has dimension $55-28=27$. All these counts
%are independent of the choice of coordinates. Hence we conclude:
%the ideal of an anticanonically embedded Severi-Brauer surface
%is generated by 27 linearly independent quadrics.

\section{The Lie algebra of a Severi-Brauer surface}\label{se:lie}

In this section we start preparing the ground for the algorithm
that establishes whether a given Severi-Brauer surface has a rational 
parametrization. 
We let $S_0$ be the Severi-Brauer surface given by
the standard embedding of $\p^2$~(\ref{SBVpoint}).
Furthermore, $S$ will be an arbitrary Severi-Brauer surface
anticanonically embedded into $\p^9$. 
By Theorem~\ref{thm:S} $S$ and $S_0$ are isomorphic over $\Q$
exactly if there is a matrix $M$ such that
\begin{equation}\label{M}
  M\in{\GL}_{10}(\Q)\textrm{ and }
  p\mapsto Mp \textrm{ is a bijection from } S_0 \textrm{ to } S.
\end{equation}
Finding an isomorphism of a given surface $S$ with $S_0$ and hence
a parametrization of $S$ therefore means finding $M$ such that~(\ref{M}) holds.

For the moment we work over an arbitrary field $F\subset\C$.
Let the anticanonically embedded Severi-Brauer surface $S$ have 
a point over $F$ and hence be isomorphic to $\p^2(F)$.
Then its automorphism group $\Aut(S)$ is isomorphic to 
the automorphism group of $\p^2(F)$, which is $\PGL_3(F)$ 
(cf.~\cite{har92}).
In particular
$S$ admits only linear automorphisms, so $\Aut(S)$ consists
of all $g\in \PGL_{10}(F)$ such that $gp\in S$ for all $p\in S$.
We recall that an anticanonically embedded Severi-Brauer surface $S$ 
is given by $27$ independent quadrics. In 
other words, there are $27$ symmetric matrices $A_i$ such that
$p\in S$ if and only if $p^TA_ip =0$ for $1\leq i\leq 27$.
If we denote 
%$$G(S,F) = \{  g\in {\GL}_{10}(F) \mid \text{ there are $\lambda_{ij}\in F$ 
%such that } g^TA_i g = \sum_{j=1}^{27} \lambda_{ij} A_j \},$$
$$G(S,F) = \{  g\in {\GL}_{10}(F) \mid \exists\ \lambda_{ij}\in F 
\text{ s.t. } g^TA_i g = \sum_{j=1}^{27} \lambda_{ij} A_j \},$$
then $\Aut(S)\cong G(S,F)/Z$, where $Z$ is a subgroup of $G(S,F)$
consisting of all scalar matrices. However, we rather work with the
group $G(S,F)$, since it is conveniently given by $10\times 10$-matrices.
Also we set
%$$L(S,F) = \{ X\in \mf{gl}_{10}(F) \mid \text{ there are $\lambda_{ij}\in F$ 
%such that } X^TA_i +A_iX = \sum_{j=1}^{27} \lambda_{ij} A_j \}.$$ 
$$L(S,F) = \{ X\in \mf{gl}_{10}(F) \mid \exists\ \lambda_{ij}\in F 
\text{ s.t. } X^TA_i +A_iX = \sum_{j=1}^{27} \lambda_{ij} A_j \}.$$ 

\begin{lemma}\label{realLieAlg}
The group $G(S,\R)$ is a Lie group and its Lie algebra is $L(S,\R)$.
\end{lemma}

\begin{proof}
%For an overview of the theory of Lie groups we refer to~\cite{knapp02}.
To see that $G(S,\R)$ is a Lie group as described in~\cite{knapp02}, 
we embed it into $\GL_{37}(\R)$.
For $g\in G(S,\R)$ let $\lambda_{ij}\in\R$ be such that 
$g^T A_i g = \sum_{j=1}^{27}\lambda_{ij}A_j$,
then $g$ is mapped to $\mathrm{diag}(g,\Lambda)$ 
where $\Lambda=(\lambda_{ij})$ is the matrix containing $\lambda_{ij}$.

We use the characterization of the Lie algebra of $G(S,\R)$ as the
set of all $X\in M_{10}(\R)$ such that $\exp( tX )\in G(S,\R)$
for all $t\in \R$ (cf., \cite{knapp02}). Let $X$ have this property.
Then there are $\nu_{ij}(t)$ such that $(\exp( tX) )^TA_i
\exp(tX) = \sum_j \nu_{ij}(t) A_j$ for $t\in \R$. We differentiate
this equation with respect to $t$ and set $t=0$. This yields $X\in 
L(S,\R)$. \par
Suppose on the other hand that $X \in L(S,\R)$, and let $\lambda_{ij}
\in\R$ be such that $X^TA_i +A_iX = \sum_{j=1}^{27} \lambda_{ij} A_j$.
For $s\geq 0$ and $t\in \R$ set
$$R_s = \sum_{r=0}^s \frac{(tX^T)^r}{r!} A_i \frac{(tX)^{s-r}}{(s-r)!}.$$
Then $(\exp tX)^T A_i (\exp tX) = \sum_{s\geq 0} R_s$. A small 
calculation shows
$$R_{s+1}= \frac{t}{s+1}(X^TR_s+R_sX).$$
Let $\Lambda=(\lambda_{ij})$ be the matrix containing the $\lambda_{ij}$.
Then by induction we get $R_s = (1/s!) \sum_j t^s\Lambda^s(i,j)A_j$.
Hence $(\exp tX)^T A_i (\exp tX) = \sum_j \exp( t \Lambda)(i,j) A_j$.
It follows that $X$ lies in the Lie algebra of $G(S,\R)$.
\end{proof}

From the discussion before the Lemma it follows that $L(S,\R)/K$ with 
$K$ the subalgebra consisting of the scalar matrices in $L(S,\R)$, is the
Lie algebra of the group of automorphisms of $S$. 
%For our purpose it is more convenient to work with whole $L(S,\R)$

\begin{remark}
There is also an alternative way for finding the Lie algebra of $G(S,F)$
corresponding to $S$ by understanding $G(S,F)$ as an algebraic group.
For an overview of the theory of algebraic groups we refer to~\cite{hum81}.
To describe the group $G(S,F)$ by polynomial functions, we embed it into 
$\GL_{37}(F)$ as in the proof of Lemma~\ref{realLieAlg}.
%Then $G$ consists of block diagonal matrices $\mathrm{diag}(g,\Lambda)$ 
%where $\Lambda=(\lambda_{ij})$ is the matrix containing the $\lambda_{ij}$.
The ideal defining the image of $G(S,F)$ under this embedding 
is generated by two types of polynomials:
\begin{align*}
  f^i_{rs}(\mathbf{t}) &= (g^TA_i g - \sum_{j=1}^{27} \lambda_{i,j} A_j)_{rs} 
      = \sum_{k,l=1}^{10}(A_i)_{kl}t_{kr}t_{ls} - 
        \sum_{j=1}^{27} t_{10+i,10+j}(A_j)_{rs} \\
       &\qquad\qquad \textrm{ for } r,s=1,\ldots,10,\ i=1,\ldots,27, \textrm{ and }\\
   g_{rs}(\mathbf{t}) &= t_{rs}
    \quad\textrm{ if } r=1,\ldots,10, s=11,\ldots,37 \textrm{ or } 
                       r=11,\ldots,37, s=1,\ldots,10.
\end{align*}
Then the Lie algebra $L$ consists of those 
$\delta = (\delta_{kl})\in\mf{gl}_{37}(F)$ that satisfy
$\delta(f^i_{rs}) = 
  \sum_{k,l=1}^{37}\delta_{kl}\frac{\partial}{\partial t_{kl}}f^i_{rs}|_{t=e} = 0$ 
(i.e.~partial derivatives followed by the evaluation at the identity $e$)
and similarly $\delta(g_{rs}) = 0$ for all relevant $r,s,i$.
This gives us conditions for $\delta$:
$\bar\delta^TA_i + A_i\bar\delta = \sum_{j=1}^{27} \delta_{10+i,10+j} A_j$,
where $\bar\delta = (\delta_{kl})_{k,l=1}^{10}$ is the left upper block of $\delta$,
and $\delta_{rs}=0$ for $r=1,\ldots,10, s=11,\ldots,37$ or 
$r=11,\ldots,37, s=1,\ldots,10$.
This is isomorphic to the algebra $L(S,F)$.
\end{remark}

\begin{theorem}\label{LieAlgL0}
Let $S_0\subset\mathbb{P}^9(\Q)$ be the Severi-Brauer surface given
by the standard embedding~(\ref{SBVpoint}). Let $L_0 = L(S_0,\Q)$.
Then $L_0$ is isomorphic to $\mf{gl}_3(\Q)$ and the natural
10-dimensional $L_0$-module is irreducible.
\end{theorem}

\begin{proof}
We first prove the statement of the theorem over $\R$, 
and afterwards we revert back to $\Q$.

As above let $Z$ denote the subgroup of $G(S_0,\R)$ consisting of the
scalar matrices. Let $K$ denote the Lie algebra of $Z$. If we view $K$ as a
subalgebra of $L(S_0,\R)$ then $K$ coincides with the scalar matrices
in $L(S_0,\R)$. Now the Lie algebra of $G(S_0,\R)/Z$ is isomorphic to 
$L(S_0,\R)/K$.  However, since $G(S_0,\R)/Z$ is isomorphic to $\PGL_3(\R)$
also the Lie algebra of $G(S_0,\R)/Z$ is isomorphic to the Lie algebra
of $\PGL_3(\R)$, which is isomorphic to $\mf{sl}_3(\R)$. It follows that
$L(S_0,\R)/K$ is isomorphic to $\mf{sl}_3(\R)$. Since $K$ is 1-dimensional
we get that $L(S_0,\R)$ is isomorphic to $\mf{gl}_3(\R)$.

Let $\{v_0,v_1,v_2\}$ be the standard basis of $V=\R^3$. Let
$W=\Sym^3(V)$ with the basis
$\{v_0^3$, $v_1^3$, $v_2^3$,
$3v_0^2v_1$, $3v_1^2v_2$, $3v_2^2v_0$,
$3v_0v_1^2$, $3v_1v_2^2$, $3v_2v_0^2$,
$6v_0v_1v_2\}$. Let $\varphi:V\to W$ be defined by
$\varphi(v) = v^3$. We write the coordinates of an element of
$W$ with respect to the basis above. Then the image of the induced
map $\varphi: \p(V)\to\p(W)$ is exactly $S_0$, see (\ref{SBVpoint}).

Write $H=\GL_3(\R)$. Then $H$ acts on $W$ by $h\cdot uvw = (hu)(hv)(hw)$,
for $u,v,w\in V$. By writing the matrix of elements of $H$ with respect to
the basis above we get a representation $\rho : H\to \GL_{10}(\R)$. We have
$h\cdot \varphi(v) = \varphi(h\cdot v)$, and hence $\varphi(V)$ is fixed
under the action of $H$ on $W$. We have further $S_0 = \varphi(\p(V))$,
therefore $\rho(H)\subseteq\Aut(\varphi(V))=G(S_0,\R)$.

Let $h=(h_{ij})_{i,j=0}^2\in H$.
Since $(\rho(h))_{ij} = h_{ij}^3$ for $i,j = 0,1,2$, we see that
$\rho:H\to G(S_0,\R)$ is injective.
Hence by differentiating $\rho$ we get an injective morphism
of Lie algebras $d\rho : \mf{gl}_3(\R)\to L(S_0,\R)$.
Since the dimensions of these Lie algebras are equal,
this is an isomorphism.
Therefore the natural $L(S_0,\R)$-module is isomorphic to $\Sym^3(V)$ and
hence irreducible.

Now we can prove the statement of the theorem for $\Q$.
Let $e_{ij}$ denote the element of $\mf{gl}_3(\R)$ that has 
a $1$ on position $(i,j)$ and zeros elsewhere. 
Since the module afforded by the representation
$d\rho : \mf{gl}_3(\R)\to \mf{gl}_{10}(\R)$ is $\Sym^3(V)$,
it follows then that the matrix $d\rho(e_{ij})$ has integer entries
as well. 
Since $L_0= d\rho(\mf{gl}_3(\R))\cap \mf{gl}_{10}(\Q)$, 
we see that $d\rho(e_{ij})\in L_0$.
Hence $L_0\cong\mf{gl}_3(\Q)$.

Finally, the natural $L_0$-module is irreducible over $\Q$ since it is
irreducible over $\R$.
\end{proof}

Now we revert back to the situation described at the beginning 
of the section, i.e.~we let $S_0,S\subset\p^9(\Q)$
be two Severi-Brauer surfaces, each given by 27 linearly independent
quadrics $x^TA_i^0x$ resp. $x^TA_ix$, $i=1,\dots 27$,
where $x=(x_0\dots x_9)^T$ is a column vector. 
The embedding of $S_0$ in $\p^9(\Q)$ is also given by
the standard embedding~(\ref{SBVpoint}).
In the sequel we work with the Lie algebras
$L_0 = L(S_0,\Q)$ and $L = L(S,\Q)$.

\begin{proposition}\label{lieAlgisom}
Suppose that $S_0$ and $S$ are isomorphic over $\Q$ and let $M$ be 
the matrix as in~(\ref{M}).
Then the map $X\mapsto MXM^{-1}$ is a Lie algebra
isomorphism from $L_0$ to $L$.
\end{proposition}

\begin{proof}
We have that $p\in S_0$ if and only if $Mp\in S$. In other words
$p^TA_i^0p =0$ for all $i$ if and only if $p^TM^TA_iMp=0$ for all $i$.
Hence the matrices $M^TA_iM$ also describe $S_0$. So there are $\eta_{ij}
\in \Q$ such that $M^TA_iM = \sum_{j=1}^{27} \eta_{ij} A_j^0$.
In the same way there are $\theta_{ij}\in \Q$ such that 
$(M^{-1})^TA_i^0M^{-1} = \sum_{j=1}^{27} \theta_{ij} A_j$.\par
Now let $X\in L_0$. Then we have to show that $MXM^{-1}\in L$. 
Using the above we calculate
\begin{align*} 
(MXM^{-1})^TA_i +A_iMXM^{-1} &= \sum_{j=1}^{27} \eta_{ij} 
(M^{-1})^T(X^TA_j^0+A_j^0X)M^{-1} \\
&= \sum_{j,k=1}^{27} \eta_{ij}\lambda_{jk} (M^{-1})^TA_k^0 M^{-1}\\
&= \sum_{j,k,l=1}^{27} \eta_{ij}\lambda_{jk}\theta_{kl} A_l.
\end{align*} 
The desired conclusion follows.
\end{proof}

We decompose $L_0$ and $L$ as direct sums of ideals as follows
$$ L_0 = \left<I_{10}\right> \oplus \{x\in L_0 \mid \Tr(x)=0\}, $$
$$ L = \left<I_{10}\right> \oplus \{x\in L \mid \Tr(x)=0\}. $$
We call the ideal $\{x\in L_0 \mid \Tr(x)=0\}$
the traceless part of $L_0$, and similarly for $L$. 
Then by Theorem~\ref{LieAlgL0} the traceless part of $L_0$ is
isomorphic to $\mf{sl}_3(\Q)$. So if $L_0\cong L$ then 
the traceless part of $L$ is isomorphic to $\mf{sl}_3(\Q)$ as well.

In the next sections we describe an algorithm that given a Lie algebra $N$
such that $N\otimes\C\cong\mf{sl}_3(\C)$ decides whether $N$ is isomorphic to
$\mf{sl}_3(\Q)$ and in the affirmative case explicitly constructs an isomorphism.
Therefore we can decide whether $L_0$ and $L$ are isomorphic. If this is not the case,
we conclude that $S_0$ and $S$ are not isomorphic over $\Q$, and we are done.
So in the remainder of this section we assume that $L_0$ and $L$ are isomorphic.
Set $K=\mf{sl}_3(\Q)$, then we construct isomorphisms of $K$ with the
traceless part of $L_0$ and $L$ respectively. This yields injective
homomorphisms $\varphi_0:K\to L_0$ and $\varphi:K\to L$. We note
that these maps are representations of $K$.

Let $H\subset K$ be a fixed Cartan subalgebra with basis $h_1,h_2$ which are
part of a Chevalley basis of $K$. Let $\tau$ 
be a fixed automorphism of $K$, such that $\tau(h_1)=h_2$
and $\tau(h_2)=h_1$ (such an automorphism exists by \cite{gra6}, \S 5.11).

\begin{theorem}
  The representation $\varphi_0$ of $K$ is either isomorphic to $\varphi$
  or to $\varphi\circ\tau$.
  Let $V_0$ be the 10-dimensional $K$-module corresponding to $\varphi_0$.
  If $\varphi_0$ is isomorphic to $\varphi$ then we let $V$ be 
  the 10-dimensional $K$-module corresponding to $\varphi$, otherwise
  we let $V$ be the $K$-module corresponding to $\varphi\circ\tau$.
  Let $f:V_0\to V$ be an isomorphism of $K$-modules.
  Then $f$ modulo scalar multiplication is also an isomorphism from $S_0$ to $S$.
\end{theorem}

\begin{proof}
First note that the $K$-modules $V_0$ and $V$ are irreducible. 
For $V_0$ this follows from Theorem~\ref{LieAlgL0}. However, since $L_0$ and
$L$ are isomorphic, the same holds for $V$ (cf.~Proposition~\ref{lieAlgisom}).

There are exactly two irreducible $K$-modules of dimension $10$. 
We represent a weight $\lambda \in H^*$ by the tuple $(\lambda(h_1),
\lambda(h_2))$. Then the two irreducible $K$-modules of dimension
$10$ have highest weights $(3,0)$ and $(0,3)$ respectively.
By composing $\varphi$ with $\tau$ we change the highest weight
of the corresponding module (from $(3,0)$ to $(0,3)$ or vice versa).
Therefore, after maybe composing $\varphi$ with $\tau$ we have that 
the two representations have the same highest weight, and hence 
are isomorphic.

In order to prove the last assertion of the theorem 
we may work over $\C$. We consider the
Lie groups $H_0$ and $H$ which are generated by all exponentiations of 
$\varphi_0(K\otimes \C)$, and all exponentiations of $\varphi(K\otimes\C)$ 
respectively. Then $H_0$ and $H$ are subgroups of $G(S_0,\C)$ and $G(S,\C)$
respectively (cf. Lemma~\ref{realLieAlg}).
Both Lie groups are isomorphic to
the Lie group generated by all exponentiations of the elements of the 
natural representation of $K\otimes \C$, which is $\SL_3(\C)$.
These two isomorphisms give two representations $\psi_0,\psi$
of $\SL_3(\C)$ in $V_0\otimes\C$ and $V\otimes\C$, respectively, 
and $f$ is an isomorphism between these two representations.

The orbits of $\SL_3(\C)$ in the representation $\psi_0$ are the orbits of
cubic forms under linear substitution. By Table~5.16 in \cite{Dimca:87},
there is precisely one $\GL_3(\C)$-orbit of dimension~3, namely the
orbit of triple lines; all other orbits except the zero orbit have
dimension at least~5. The orbit of triple lines is also a single
orbit under the action of $\SL_3(\C)$, hence it is the only 
$\SL_3(\C)$-orbit of dimension~3. Its image in the projective
space ${\mathbb P}^9$ is therefore equal to the surface $S_0$. 
Similarly, the image of the unique 3-dimensional orbit under $\psi$
in ${\mathbb P}^9$ is $S$. It is clear that $f$ takes the unique
3-dimensional orbit under $\psi_0$ to the unique 3-dimensional orbit
under $\psi$.
\end{proof}

\begin{remark}
Let $\varphi' : K\to L$ be a different embedding of $K$ into
$L$, and let $V'$ be the associated $K$-module, which we assume to
be isomorphic to $V_0$. Then there is another
way to see that an isomorphism $f': V_0\to V'$ will also lead to an
isomorphism $S_0\to S$. For that set $g = \varphi'^{-1}\circ \varphi$.
Then $g$ is an automorphism of $K$, hence by \cite{jac}, Theorem 4 of 
Chapter IX, we infer that either $g$ is an inner automorphism, or the 
composition of $\tau$ and an inner automorphism. By \cite{gra6} \S 8.5, 
composing a representation with an inner automorphism does not change its
highest weight. Therefore $g$ has to be inner.
First suppose that $g$ is of the form $\exp(\ad z)$, where $z\in K$ is 
such that $\ad z$ is nilpotent. This implies that $\varphi(z)$ is nilpotent
as well. Now by the proof of \cite{gra6}, Lemma 8.5.1, the map 
$\exp( \varphi(z) )$ is a module isomorphism $V\to V'$. 
Furthermore, since $V$ and $V'$ are irreducible $K$-modules, this is the only
isomorphism upto scalar multiples. Hence
$f' =  \exp( \phi(z) ) \circ f$ (upto scalar multiples). But 
$\exp( \phi(z) )$ lies in the automorphism group of $S$. It follows
that $f'$ also provides an isomorphism $S_0\to S$. If $g$ is a product 
of elements of the form $\exp( \ad z )$, then we reach the same conclusion.

Note that this construction generalizes to finding an isomorphism 
of Severi-Brauer varieties of arbitrary dimension $n$, since
the modules involved are always symmetric powers 
of the natural $\mf{sl}_{n+1}$-module, therefore irreducible.
\end{remark}

\begin{remark}
We note that by using standard techniques from the representation
theory of semisimple Lie algebras it is straightforward to
construct an isomorphism between two irreducible $K$-modules $V$ and
$W$. Let $\lambda$ be the highest weight of both modules, and 
$v_{\lambda}\in V$, $w_{\lambda}\in W$ be two corresponding highest
weight vectors. Let $y_1,y_2,y_3$ be the negative root vectors
of $K$, spanning the subalgebra $N^-$ of $K$. Compute a set of 
elements $u_k$ in the universal enveloping algebra of $N^-$ such 
that the elements $u_k\cdot v_{\lambda}$
form a basis of $V$. Then a module isomorphism maps $u_k\cdot v_{\lambda}$
to $u_k\cdot w_{\lambda}$.
\end{remark}

\section{Construction of a central simple associative algebra}\label{se:sln}

In this section $\gn$ will be a simple Lie algebra of 
dimension $8$ over $\Q$ such that $\gn\otimes\C$ is isomorphic to
$\mf{sl}_3(\C)$. The problem is to decide whether this isomorphism 
already exists over $\Q$, i.e., whether $\gn$ is isomorphic to
$\mf{sl}_3(\Q)$.\par
As a first step towards deciding this we compute a Cartan subalgebra
$H$ of $\gn$ (cf., \cite{gra6}, \S 3.2). The problem is
immediately solved if $H$ happens to be split (i.e., $\ad h$ has all 
its eigenvalues in $\Q$ for all $h\in H$). In that case
we can effectively construct an isomorphism $\gn \to \mf{sl}_3(\Q)$, 
for example by using the method of \cite{gra6}, \S 5.11. So in the
remainder of this section we suppose that $H$ is not split. \par
The second step of our method is the construction of a certain rational
representation of $\gn$. This works as follows. Let $F$ denote a number 
field containing the eigenvalues of $\ad h$ for $h\in H$. Consider the Lie 
algebra $\gn\otimes F$. This Lie algebra has a split Cartan subalgebra 
(namely $H\otimes F$). Therefore we can find an isomorphism 
$\gn\otimes F \to \mf{sl}_3(F)$. This gives us a representation
$\rho' : \gn\otimes F\to \mf{gl}(V')$, where $V'$ is a vector space over $F$
of dimension $3$. Now we view $V'$ as a vector space over $\Q$. More
precisely, let $W$ be the $\Q$-span of a basis of $V'$. Since $F$ is a 
vector space over $\Q$ we can form the tensor product $V = 
F\otimes_{\Q} W$, which is a vector space over $\Q$. There is a 
bijective $\Q$-linear map $\phi : V\to V'$ with $\phi(\alpha\otimes w)
= \alpha w$. From this we get a representation $\rho : \gn \to \mf{gl}(V)$,
where $\rho(x) = \phi^{-1} \circ \rho'(x) \circ \phi$. \par
By $\rho(\gn)^*$ we denote the associative algebra over $\Q$ generated
by $\rho(\gn)$. \par
For an associative algebra $A$ we let $A_{\mathrm{Lie}}$ be
the Lie algebra associated to $A$
(i.e., it has the same underlying vector space as $A$, and the Lie
product is formed by $[x,y]=xy-yx$).
The following lemma is immediate.

\begin{lemma}\label{le:lie2assoc}
  Set $A = \rho(\gn)^*$.
  Then $\rho : \gn \to A_{\mathrm{Lie}}$
  is an injective homomorphism of Lie algebras.
\end{lemma}

\begin{lemma}
Suppose that $\gn$ is isomorphic to $\mf{sl}_3(\Q)$. 
Then $\dim _{\Q}\rho(\gn)^* = 9$.
\end{lemma}

\begin{proof}
Since $\gn$ is isomorphic to $\mf{sl}_3(\Q)$ we have that $\gn$ has a
$\Q$-basis $x_1,\ldots,x_8$ that is a Chevalley basis (cf. \cite{hum}, 
\S 25.2). For example we can take the image of a Chevalley basis 
of  $\mf{sl}_3(\Q)$ under a $\Q$-isomorphism $\mf{sl}_3(\Q) \to
\gn$. \par
By \cite{hum} Theorem 27, $V'$ has a basis $B$ such that the matrix 
with respect to $B$ of each $\rho'(x_i)$ has integer entries. In other 
words, there exists an $X\in  \mathrm{GL}_3(F)$ such that 
$X^{-1}\rho'(x_i)X$ has integer entries, and hence lies in $\mf{sl}_3(\Q)$. 
So $I_3$ (the $3\times 3$-identity matrix) along with the $X^{-1}\rho'(x_i)X$ 
form a basis of $M_3(\Q)$ (the algebra of $3\times 3$-matrices over $\Q$). 
Hence there are $c_{ij}^k\in \Q$ and $e_{ij}\in \Q$ such that 
$$ X^{-1}\rho'(x_i)X X^{-1}\rho'(x_j)X = \sum_{k=1}^8 c_{ij}^k
X^{-1}\rho'(x_k)X + e_{ij}I_3,$$
or
$$\rho'(x_i)\rho'(x_j) = \sum_{k=1}^8 c_{ij}^k \rho'(x_k) + e_{ij}I_3.$$
Hence the $\rho'(x_i)$ along with $I_3$ span a $9$-dimensional associative
algebra over $\Q$. We finish with the observation that the $\Q$-dimension
of $\rho(\gn)^*$ is equal to the $\Q$-dimension of $\rho'(\gn)^*$
(since $\rho(\gn)^* = \phi^{-1} \rho'(\gn)^*\phi$).
\end{proof}

In the third step of the method we check whether the dimension 
of $\rho(\gn)^*$ is $9$. If not, then $\gn$ is not isomorphic to
$\mf{sl}_3(\Q)$ and we are done. If this dimension is $9$ then we
proceed. 

\begin{proposition}
Set $A=\rho(\gn)^*$, and suppose that $\dim_{\Q}(A) = 9$. Then
$A$ is a central simple algebra. Furthermore, $\gn$ is isomorphic
to $\mf{sl}_3(\Q)$ if and only if $A$ is isomorphic to $M_3(\Q)$
(the algebra of $3\times 3$-matrices over $\Q$).
\end{proposition}

\begin{proof}
Set $K = \rho(\gn)$.
Then by Lemma~\ref{le:lie2assoc} $K$ is isomorphic to $\gn$.
It follows that there is a direct sum decomposition $A_{\mathrm{Lie}} =
K\oplus Z$, where $Z$ is spanned by the identity of $A$. Now any two-sided
ideal of $A$ is an ideal of $A_{\mathrm{Lie}}$ as well. However, the only 
ideals of $A_{\mathrm{Lie}}$ are $0$, $Z$, $K$ and $A_{\mathrm{Lie}}$. 
Now $Z$ is not an ideal of $A$. So the only ideal of $A$, not equal
to $0$ or $A$, has to coincide with $K$. This ideal cannot be nilpotent
(otherwise $\gn$ would be a nilpotent Lie algebra by Engel's theorem,
cf. \cite{hum}),
hence the radical of $A$ is zero. So $A$ is the direct sum of simple
ideals. But the centre of $A$ has dimension $1$ (it is equal to the
centre of $A_{\mathrm{Lie}}$), so there is only one simple ideal in
the direct sum decomposition of $A$. We conclude that $A$ is
central simple.\par
If $A$ is isomorphic to $M_3(\Q)$, then $A_{\mathrm{Lie}}$ is isomorphic to
$\mf{sl}_3(\Q) \oplus Z$. Hence $\gn$ is isomorphic to $\mf{sl}_3(\Q)$.
For the other direction, if $\gn$ is isomorphic to $\mf{sl}_3(\Q)$, then
$\gn$ contains a split Cartan subalgebra $\widetilde{H}$. There is a 
basis of $V$ with respect to which the matrix of $\rho(h)$ is diagonal
for all $h\in \widetilde{H}$ (cf. \cite{gra6}, \S 8.1). 
Adding the identity we see
that $A$ contains a $3$-dimensional split torus (i.e, a commutative 
diagonalisable subalgebra). Now because $A$ is central simple, it is 
isomorphic to $M_r(D)$, where $D$ is a division ring over $\Q$. But such an 
algebra only contains a $3$-dimensional split torus if $r=3$ and $D=\Q$.
\end{proof}

%\begin{remark}
%  How unique is the embedding $\gn \to A_{\mathrm{Lie}}$?
%  The question is answered in~\cite{jac}, \S 10, Theorem~10.
%  If $\gn$ is of type $A_I$, then there are precisely two such embeddings. 
%  If $\gn$ is of type $A_{II}$, then there is no such embedding.
%\end{remark}

\begin{remark}
  A question arises: is an embedding $\gn \to A_{\mathrm{Lie}}$ unique? 
  This is answered in~\cite{jac}, Chapter~10, Theorem~10, namely, 
  there are either precisely two such embeddings or there are none.
\end{remark}

Now let $A$ be as in the lemma, and suppose that $A$ is isomorphic 
to $M_3(\Q)$. Let $\tau : M_3(\Q) \to A$ be an isomorphism. Then we take
the associated Lie algebras of these associative algebras, and restrict 
$\tau$ to the subalgebra $\mf{sl}_3(\Q)\subset (M_3(\Q))_{\mathrm{Lie}}$.
We obtain an embedding $\tau : \mf{sl}_3(\Q) \to A_{\mathrm{Lie}}$. The 
image of this map is the unique semisimple subalgebra of 
$A_{\mathrm{Lie}}$, i.e., $\rho(\gn)$. Hence after composing with 
$\rho^{-1}$ we obtain an isomorphism $\mf{sl}_3(\Q) \to \gn$.

\section{Associative central simple algebra}\label{se:assoc}

In this section we start with a central simple algebra $A$ of degree 3
and we want to decide whether $A$ is isomorphic to $M_3(\Q)$. In the
affirmative case we also want to construct an isomorphism.

Since the degree of the algebra $A$ is prime, 
by Wedderburn's structure theorem there are only two possibilities: 
either $A \cong M_3(\Q)$ or $A$ is a division algebra.
So if we by chance find a zero divisor in $A$, then we can already conclude, 
that $A\cong M_3(\Q)$. Using the zero divisor we can even
explicitly construct an isomorphism.

Let $a\in A$ be a zero divisor, so $A\cong M_3(\Q)$. 
We can find a 3-dimensional left ideal. Namely,
the vector space endomorphism $\rho_a$ of $A$: $x\mapsto xa$ 
has a nontrivial kernel (because $a$ is a zero divisor) 
and a nontrivial image (because $1.a \neq 0$). 
Both $\Ker \rho_a$ and $\im \rho_a$ are left ideals in $A$. 
Since any minimal left ideal of $M_3(\Q)$ is of dimension 3
and any left ideal is a direct sum of minimal left ideals, 
we have either $\dim(\Ker \rho_a)=3$ or $\dim(\im \rho_a)=3$.

Let $\mathcal{L}$ be a 3-dimensional left ideal in $A$.
Let $B = (b_1\ b_2\ b_3)^T$ be the column vector containing a basis 
of $\mathcal{L}$. Let $\varphi:A\to M_3(\Q)$ be the map that assigns 
to $x\in A$ the transpose of the matrix of its left action 
on $\mathcal{L}$ w.r.t.~$B$, 
so $\varphi(x) = X$, if $X^T B = (xb_1\ xb_2\ xb_3)^T$.

\begin{theorem}\label{the:zd}
  $\varphi$ is an isomorphism of algebras.
\end{theorem}

\begin{proof}
  Elementary computations show that $\varphi$ is a homomorphism of algebras. 
  Further $\varphi$ is a bijection, since otherwise $\Ker\varphi\neq 0$ would
  be a nontrivial ideal in $A$.
\end{proof}

\begin{remark}
The problematic part of the foregoing construction is the assumption
that we have a zero divisor in the algebra and so we can construct 
a 3-dimensional left ideal.
There is an ongoing research on this topic
(\cite{stoll}, the case of algebras of degree 2 was solved in~\cite{gabor2}).
First, a maximal order in the algebra is constructed 
so the structure constants are integral (cf.~\cite{gabor1}). 
Afterwards, the basis
of the order is changed to reduce the size of the structure constants
substantially. This makes finding a zero divisor possible.
Here we follow another approach,
which is similar to the one in~\cite{acc96}.
Though, the method of~\cite{acc96} does not provide
an explicit isomorphism.
We also note that by finding a maximal order we can decide 
whether a given central simple algebra is isomorphic
to the full matrix algebra (see~\cite{lajos}). 
However, this does not yield an explicit isomorphism.
\end{remark}
For finding an isomorphism of 
the algebra $A$ and $M_3(\Q)$ we first express $A$ as a cyclic algebra 
of degree 3, i.e. we find two elements $a, b$ 
of $A$ so that the algebra generated by them satisfies the following
\begin{definition}\label{def:cyclic}
  $A$ is a {\em cyclic algebra} of degree 3 if there are
  elements $a,b \in A$ so that
  $$ A = \left<1, a, a^2, 
                            \ b, ab, a^2b, 
                            \ b^2, ab^2, a^2b^2\right>_\Q, $$
  and the multiplication rules satisfy the following conditions:
  \begin{itemize}
    \item[(i)]   E = $\left<1, a, a^2\right>_\Q$ is a Galois extension of $\Q$
      with $\Gal(E|\Q) = \left<\sigma\right>$,
    \item[(ii)]  $ba = \sigma(a)b$,
    \item[(iii)] $b^3 = \beta 1$, where $\beta\in\Q^*$.
  \end{itemize}
  In this case we write $A=(E,\sigma,\beta)$.
\end{definition}

Any cyclic algebra is central simple. 
For algebras of degree 3 also the other direction holds: 
any central division algebra of degree 3 is cyclic
(a result due to Wedderburn) 
and likewise the split algebra $M_3(\Q)$ 
($(E,\sigma,1)$ is split for any $E$).
Hence in any case there exists an isomorphism of $A$ and a cyclic algebra.

In the construction of an isomorphism
the most difficult step is to find a cyclic element, 
i.e.~an element $a \in A$ 
such that the minimal polynomial $m_a(\lambda) \in \Q[\lambda]$ is
irreducible of degree 3 and the splitting field of $m_a$ has dimension
3 over $\Q$. For finding a cyclic element
we follow~\cite{jac96} 
(Lemma~2.9.8 and Theorem 2.9.17, pp.~68-70).
Although the construction is originally designed for division algebras,
it can be partially carried out also in the case $A\cong M_3(\Q)$.
In fact, the only complication which could pop up is, that we hit 
a noninvertible element and so can not continue with the original
computation. But in such case we have found a zero divisor and we
can construct the isomorphism $A\iso M_3(\Q)$ as in 
Theorem~\ref{the:zd}.

We will need the following properties of the elements 
in the matrix algebra.

\begin{lemma}\label{lem:conj}
  For a field $F$ let $A \cong M_3(F)$ and 
  let $a\in A$, $a$ not central. 
  Let the minimal polynomial $m_a(\lambda)\in F[\lambda]$ of $a$
  be irreducible over $F$.
  Then $\deg m_a(\lambda) = 3$ and
  every $b\in A$, such that $m_a(b)=0$ is a conjugate of $a$.
\end{lemma}

\begin{proof}
  $a\notin Z(A)$ implies $\deg m_a(\lambda) > 1$.
  The case $\deg m_a(\lambda) = 2$ is not possible. 
  For, let $\deg m_a(\lambda) = 2$, $m_a(\lambda)$ irreducible. 
  Then the characteristic polynomial $\chi_a(\lambda)$ of $a$ is
  $\chi_a(\lambda) = m_a(\lambda)l(\lambda)$ with $l(\lambda)$ linear.
  The factor $l(\lambda)$ of the characteristic polynomial of $a$ 
  is irreducible, therefore it divides the minimal polynomial
  $m_a(\lambda)$, a contradiction.
  Since $m_a(\lambda)$ is irreducible over $F$ and $m_a(b)=0$, 
  it is also the minimal polynomial of $b$. 
  Then $I\lambda-a$ and $I\lambda-b$ have the same invariant factors 
  and hence $a$ and $b$ are conjugate (cf.~\cite{smi65}).
\end{proof}

Now we want to construct a cyclic element $a\in A$.
We start by picking a noncentral element $x\in A$. 
If the minimal polynomial $m_x(\lambda)$ of $x$ is reducible 
over $\Q$, we can construct a zero divisor and hence an isomorphism
$A \iso M_3(\Q)$. 
Indeed, let $m_x(\lambda) = m_1(\lambda)m_2(\lambda)$, 
then $m_1(x), m_2(x) \neq 0$ are zero divisors. 
On the other hand, if $m_x(\lambda)$ is irreducible, 
then $\deg m_x(\lambda) = 3$. 
For $A\cong M_3(\Q)$ it follows from Lemma~\ref{lem:conj}.
If $A$ is a central division algebra, then 
$m_x(\lambda)$ is the minimal polynomial of $x$ also after
extending the field of coefficients. If $F$ is a splitting field
of $A$, so $A\otimes_\Q F\cong M_3(F)$, then again by 
Lemma~\ref{lem:conj} $\deg(m_x)=3$. 
Therefore in our construction we can assume that 
the minimal polynomial of a randomly chosen element 
$x\in A$ is irreducible and has degree 3.

Recall that $[a,b] = ab-ba$.

\begin{lemma}\label{lem:auxConstr1}
  Let $A$ be a central simple algebra over $\Q$.
  Let $x\in A$ be a noncyclic element such that 
  its minimal polynomial $m_x(\lambda)$ is cubic and irreducible.
  Then there exists $t\in A$ such that $[[t,x],x]\neq 0$.
\end{lemma}

\begin{proof}
  Let $[[t,x],x]=0$ for all $t\in A$.
  For the inner derivation $i_x:y\mapsto[y,x]$ we have $i_x^2=0$,
  where $i_x^2:y\mapsto [[y,x],x] = yx^2 - 2xyx + x^2y$.
  But the vector space endomorphisms $y\mapsto yx^2$, $y\mapsto xyx$ and 
  $y\mapsto x^2y$ are linearly independent (see~\cite{jacBA2}, the
  proof of Theorem 4.6, p. 218), a contradiction.
\end{proof}

For given $x$, $[[t,x],x]\neq 0$ is a linear condition on $t$,
so it is easy to find a desired $t\in A$.

\begin{lemma}\label{lem:auxConstr2}
  Let $x,t\in A$ be as in Lemma~\ref{lem:auxConstr1}
  and let $[t,x]$ be invertible.
  Let $y=[t,x]x[t,x]^{-1}$ and $z=[x,y]$.
  Then $z\neq 0$, $z\notin\Q$, $z^3\in\Q$.
\end{lemma}

\begin{proof}
  The assertion follows from~\cite{jac96} \S 2.9., Lemma~2.9.8.
  There it is proven for central division algebras, but the proof
  for the case $A\cong M_3(\Q)$ is almost a word-by-word copy of the
  original one. The proof is technical and therefore we omit it here.
\end{proof}

\begin{theorem}
  Let $A$ be a central simple algebra over $\Q$ of degree 3. 
  We can either find a zero divisor in $A$ or 
  we can find elements $a,b\in A$ generating the algebra $A$ 
  as a cyclic algebra as described in Definition~\ref{def:cyclic}.
\end{theorem}

\begin{proof}
  Let $x\in A$ be an arbitrary element.
  Then either $x$ fulfils the assumptions of Lemma~\ref{lem:auxConstr1},
  or $x$ is a zero divisor and we are done, or $x$ is a cyclic element.
  In this case ($x$ cyclic) we set $a=x$ and denote $E=\left<1,a,a^2\right>_\Q$,
  the maximal subfield generated by $a$. 
  Let $\sigma$ denote a generator of $\Gal(E|\Q)$.
  By factoring the minimal polynomial 
  $m_a(\lambda)$ over $E$ we find $\sigma(a)$ and afterwards an element $b\in A$ 
  s.t.~$bab^{-1} = \sigma(a)$, i.e.~$b\notin E$. 
  By Lemma~\ref{lem:conj} such $b$ exists if $A\cong M_3(\Q)$, 
  for $A$ a division algebra we refer to~\cite{jac96}.
  We claim that $b$ satisfies also the rule~(iii) of Definition~\ref{def:cyclic}.
  Indeed, by multiple conjugation we get $b^3ab^{-3}=a$, so $b^3\in Z_A(E)=E$.
  If $b^3\notin\Q$, then $b^3$ would generate $E$. But in such case $b^3b = bb^3$
  would imply $b\in Z_A(E)=E$, a contradiction.

  Now we will construct elements $a$, $b\in A$ and prove that they match
  De\-fi\-ni\-tion~\ref{def:cyclic} in case $x$ is noncyclic and not a zero divisor. 
  Note, that we may assume that the
  minimal polynomials of $a$ and $b$ are irreducible, otherwise 
  we would have found a zero divisor.
  First we construct $z$ 
  as in Lemma~\ref{lem:auxConstr2} (unless $[t,x]$ is a zero divisor, 
  in which case we are done) and denote it by $z_x$. We set $b=z_x$.
  Then we carry out the construction of the Lemma~\ref{lem:auxConstr2} again,
  using the element $b$ as $x$ to get the corresponding $z$ which we denote 
  here $z_b$. Set $a=bz_b$. By~\cite{jac96} (proof of Theorem~2.9.17, p.69) 
  $a$ is a cyclic element not in $\Q$.
  If we define $\sigma$ on $E$: $\sigma(a) = bab^{-1}$,
  a technical computation shows that $\sigma(a)a = a\sigma(a)$, 
  i.e.~$\sigma(a)\in E$, so $\sigma$ is a well-defined 
  automorphism of $E$.
  The elements $a, b$ satisfy the properties of Definition~\ref{def:cyclic} 
  and hence generate a cyclic algebra $(E,\sigma,\beta)$.
\end{proof}

If we were not lucky enough so far to hit a zero divisor in the algebra $A$
and to use it for constructing an isomorphism $A \iso M_3(\Q)$
as described in the beginning of the section, then now we have
an isomorphism of $A$ and a cyclic algebra $(E,\sigma,\beta)$.
For any cyclic algebra of degree 3 
we can already decide whether it is split or 
a division algebra over $\Q$, as can be seen in the following

\begin{theorem}
  $(E,\sigma,\beta) \cong M_3(\Q)$ if and only if there exists $x\in E$
  such that
  \begin{eqnarray}\label{eqn:norm}
    x\sigma(x)\sigma^2(x) = \frac{1}{\beta}.
  \end{eqnarray}
\end{theorem}

\begin{proof}
  The algebra $(E,\sigma,\beta) \cong M_3(\Q)$ 
  exactly if we can find a 3-dimensional left ideal.
  Let $\mathcal{L}$ be such an ideal in $(E,\sigma,\beta)$ and let 
  $a_0^\prime 1 + a_1^\prime b + a_2^\prime b^2 \in \mathcal{L}$, 
  $a_i^\prime \in E$.
  At least one of $a_i^\prime$'s is nonzero, let it be $a_0^\prime$ 
  (the other cases are treated in the same way).
  Then $1 + a_1 b + a_2 b^2 \in \mathcal{L}$.
  After multiplying from the left by elements from $E = \left<1,a,a^2\right>_\Q$ 
  we see that 
  $$ \mathcal{L} = \left<1 + a_1b + a_2b^2, 
                         a(1 + a_1b + a_2b^2), 
                         a^2(1 + a_1b + a_2b^2)\right>_\Q. $$
  Since $\mathcal{L}$ is a left ideal, then also
  $b(1 + a_1b + a_2b^2) = \beta\sigma(a_2) + b + \sigma(a_1)b^2 \in \mathcal{L}$.
  This holds, if $\beta\sigma(a_2) : 1 : \sigma(a_1) = 1 : a_1 : a_2$,
  which can be satisfied if and only if $a_1$ is a solution 
  to the norm equation~(\ref{eqn:norm}).
\end{proof}

In our algorithm we use an existing routine in Magma to test
whether the norm equation~(\ref{eqn:norm}) is solvable.
In the affirmative case we use a solution for constructing an isomorphism
of $(E,\sigma,\beta)$ and $M_3(\mathbb{Q})$. Namely, if $x\in E$ is a solution, 
then
$$ \mathcal{L} =\left< 1+xb+x\sigma(x)b^2, 
                     a(1+xb+x\sigma(x)b^2), 
                   a^2(1+xb+x\sigma(x)b^2) \right>_\Q$$ 
is a 3-dimensional left ideal and we can find an isomorphism of $A$
and $M_3(\Q)$ as in the Theorem~\ref{the:zd}.

\end{document}